\documentclass{article}%
\usepackage{amssymb,amsmath,amsthm}%
\usepackage{xcolor}
\providecommand{\U}[1]{\protect\rule{.1in}{.1in}}
\newcommand\E{{\mathbb {E}}}
\def\cA{{\mathcal{A}}}

\newtheorem {Theorem}{Theorem}[section]
\newtheorem {Proposition}{Proposition}[section]
\newtheorem {Corollary}{Corollary}[section]

\newtheorem{Remark}{Remark}[section]

\newcommand\N{{\mathbb{N}}}

\newcommand\beq{\begin{equation}}
\newcommand\eeq{\end{equation}}

\begin{document}

\title{Deviation and concentration inequalities for dynamical systems with subexponential decay of correlation.}

\author{C. Cuny\footnote{Christophe Christophe Cuny, Univ Brest, LMBA, UMR 6205 CNRS, 6 avenue Victor Le Gorgeu, 29238 Brest}, J. Dedecker\footnote{J\'er\^ome Dedecker, Universit\'e de Paris, CNRS, MAP5, UMR 8145,
45 rue des  Saints-P\`eres,
F-75006 Paris, France.}
and
F. Merlev\`ede \footnote{Florence Merlev\`ede,  Univ Gustave Eiffel, Univ Paris Est Cr\'eteil, LAMA, UMR 8050 CNRS,  \  F-77454 Marne-La-Vall\'ee, France.}}

\maketitle

\begin{abstract}
We obtain large and moderate deviation estimates, as well as concentration inequalities,  for  a class of nonuniformly expanding maps with stretched exponential 
decay of correlations.  In the large deviation regime, we also exhibit  examples  showing that the obtained upper bounds are essentially optimal.  
\end{abstract}

\noindent{\it  MSC2020 subject classifications: }  37A25; 37A50; 60F10. \\
\noindent{\it Keywords: } Large deviations, Moderate deviations, Nonuniformly expanding maps, Young towers.

\section{Introduction}

Let $X$ be a bounded metric space and $T \colon X \to X$ be a transformation,
preserving a Borel probability measure $\nu$.
Suppose that $\varphi \colon X \to {\mathbb R}$ is a ${\mathbb L}^1(\nu) $ observable with
$\int \varphi \, d\nu = 0$. We consider the  Birkhoff sums
\begin{equation}
  \label{eq:Sn}
  S_n(\varphi)  = \sum_{k=0}^{n-1} \varphi \circ T^k
\end{equation}
as a discrete time random process, defined on the probability space $(X, \nu)$. The Birkhoff ergodic theorem asserts that  if $T$ is ergodic then $\lim_{n \rightarrow \infty}n^{-1} S_n(\varphi) =0$  $\nu$-a.s.  Deviation inequalities aim to quantify the rate at which the Birkhoff sum goes to zero. One of the purposes of  this paper is then to derive estimates of the quantity 
\begin{equation}\label{lmdp-deviation}
\nu  \Big ( \sqrt{\frac{a_n}{n}} |S_n( \varphi ) |  \geq x  \Big )  \, , \, x >0 \, , 
\end{equation}
when $a_n \rightarrow 0$, in either the large deviation regime (so when $a_n = 1/ n$) or in the  moderate deviation regime (so when $na_n \rightarrow \infty$ with $a_n \rightarrow 0$).

In this paper, we shall consider nonuniformly expanding maps in the sense recalled in Section \ref{Sec:mixing}. It is well-known that the study of the deviation probability  \eqref{lmdp-deviation} is widely linked with the decay of correlations of some observables of the iterates of $T$ (see for instance \cite{Mel09,AFLV11}),  which is actually linked with the moments of the return time to the basis  of the induced map  (see again Section \ref{Sec:mixing} for the definition of the induced map and its  associated return times).  When the return time to the basis has an exponential moment (this is the case for instance of the dispersing billards, see \cite{Ch99}) and $\varphi$ is an H\"older observable, it has been proved  in  \cite{RY}  (see their Theorems B and C)  that the following  large/moderate deviation principles hold:  there exists a strictly convex function $I (x)$ vanishing only at $x=0$ such that 
\[
\lim_{n \rightarrow \infty}a_n  \log \nu  \Big ( \sqrt{\frac{a_n}{n}} S_n( \varphi )  \geq x  \Big ) =- I(x) \, .
\]
In the large deviation regime (i.e. $a_n =1/n$) we refer also to  Theorem 2.1 in \cite{MN08}.

Next, when  the return time to the basis has a moment of polynomial  order and $\varphi$ is either an H\"older observable or a bounded variation function, the deviation probability is still well understood both in the large deviation regime and in the moderate one. Indeed, it follows from \cite{DGM} that if the return time $R$ to the basis has a weak moment of order  $p>1$ (so when $m ( R > n ) \leq C n^{-p}$ where $m$ is the reference probability measure defined in Section \ref{Sec:mixing}), one has for any $x >0$, 
\[
\limsup_{n \rightarrow \infty} n^{\alpha p -1}  \nu  \Big (   |S_n( \varphi ) |  \geq x n^{\alpha } \Big )  \leq C x^{-p}\, ,
\]
for any $\alpha >1/2$ such that $1/p \leq \alpha \leq 1$. Moreover it has been proved in \cite{DGM} that this upper bound cannot be essentially improved. 

However, when the  return time $R$ to the basis has subexponential (also called stretched exponential) moment of order $\gamma \in ]0,1[ $ (so when $m ( R > n ) \leq C {\rm e}^{- \delta n^\gamma}$ for some $\delta >0$), the situation is not so well understood.  The recent paper \cite{AF} considers this case, but when the dynamical system can be modeled by a Young tower, the  obtained bounds turn out to be suboptimal (see our Remark \ref{RemarkonAF}).  In Section \ref{DIMDP}, we improve the estimates of the quantity \eqref{lmdp-deviation} both in the large deviation regime and in the moderate one, when the return time  to the basis has stretched moment of order $\gamma \in ]0,1[ $ and when the observable $\varphi$ is either  H\"older continuous or with bounded variations. For instance, in case when $T$ is the Viana map  as introduced in \cite{Via97}, using that its associated return time to the basis exhibits stretched exponential moment of order at least $1/2$ (so that $\gamma =1/2$) as proved in \cite{Gou06},  it follows from our Corollary \ref{ineexpo}  that, for any H\"older observable $\varphi$ there exist  positive constants $c_1$ and $c_2$ such that for  any $x >0$,  
\[ \limsup_{n \rightarrow \infty} \frac{1}{ \sqrt{n}}\log \nu  \Big (\frac{S_n( \varphi )}{n}  \geq x  \Big ) \leq - c_1 x^{1/2}  \,  \text{ and }  \, \limsup_{n \rightarrow \infty} a_{n}\log \nu \Big (\sqrt{\frac{a_{n}}{n}} S_n( \varphi )  \geq x  \Big ) \leq - c_2 x^2 \, , 
\]
where 
$a_{n}\rightarrow 0$ and $a_{n}n^{1/3}\rightarrow \infty $.  

Our strategy of proof will be first to estimate the $\tau$-mixing coefficients associated with the observable of the nonuniformly expanding map (see our Section \ref{Sec:mixing} for the definition of these coefficients) and then to apply previous known results for $\tau$-mixing sequences. Moreover, considering  the family of interval maps $T_{\gamma}$, $\gamma \in \,  ]0,1]$,  introduced in \cite[Appendix A]{CDKM2} and which are such that the return time to the basis has stretched exponential moment (see also our Section \ref{sec:example}) we will show that the large deviation upper bounds given in this paper are essentially optimal. Concerning the moderate deviation regime, we will also give in the Corollary 
\ref{MDP} a moderate deviation principle, which implies \eqref{mdp1} with a restriction (depending on $\gamma$) on the possible range of the sequence $a_n$. 

We shall also be interested in proving concentration inequalities in the spirit of those obtained in \cite{CG12,GM14} but when the  return time  to the basis has stretched exponential moment of order $\gamma \in ]0,1[ $ (the case of exponential moment and  strong polynomial moments is handled   in \cite{CG12}; the case of weak polynomial moments is considered in \cite{GM14}).  Section \ref{Sec:CI} is devoted to the statement of these  concentration inequalities in case of stretched exponential moment of the return time. All the proofs are postponed to Section \ref{Sec:proofs}. 

\section{Mixing properties of nonuniformly expanding maps}  \label{Sec:mixing}

\setcounter{equation}{0}

Let $(X,d)$ be a complete bounded separable metric space with the Borel $\sigma$-algebra. Let us introduce the class of dynamical systems that we consider in this paper. 
Suppose that $T \colon X \to X$ is a measurable transformation which
admits an inducing scheme consisting of:
\begin{itemize}
  \item a closed subset $Y$ of $X$ with a \emph{reference} probability measure $m$ on $Y$;
  \item a finite or countable partition $\Gamma = \cup_{\alpha \in E} \Gamma_{\alpha}$  of $Y$ (up to a zero measure set)
    with $m(a) > 0$ for all $a \in \Gamma$;
  \item an integrable \emph{return time} function $R  \colon Y \to \{1,2,\ldots\}$ which is constant
    on each $a \in \Gamma$ with value $R(a)$ and
    $T^{R (a)}(y) \in Y$ for all $y \in a$, $a \in \Gamma$.  We require in addition that ${\rm gcd} \{ R(y) , y \in Y \} = 1 $. 
\end{itemize}

Let $F \colon Y \to Y$, $F(y) = T^{R(y)}(y)$ be the induced map.
We assume that there are constants $\lambda > 1$ and  $K > 0$ 
such that for each $a \in \Gamma$ and all $x,y \in a$:
\begin{itemize}
  \item $F$ restricts to a (measure-theoretic) bijection from $a$ to $Y$;
  \item $d(F(x), F(y)) \geq \lambda d(x,y)$;
  \item $d(T^k(x), T^k(y)) \leq K d(F(x), F(y))$
    for all $0 \leq k \leq R(a)$;
  \item the inverse Jacobian $\zeta_a = \frac{dm}{dm \circ F}$ of
    the restriction $F \colon a \to Y$ satisfies
    \[
      \bigl| \log |\zeta_a(x)| - \log |\zeta_a(y)| \bigr|
      \leq K d(F(x), F(y))
      .
    \]
\end{itemize}

In addition to the standard assumptions above, we rely on
non-pathological coding of orbits under $F$ allowed by the elements
of $\Gamma$.
Let $\cA$ be the set of all finite words in the alphabet $\Gamma$
and $Y_w= \cap_{k=0}^n F^{-k}(a_k)$ for $w=a_0\cdots a_n \in \cA$.
We require that
\beq \label{nonpathologicalH}
  m(Y_w)=m(\bar Y_w)
  \qquad \text{for every } w \in \cA
  .
\eeq

We say that the map $T$ as above is nonuniformly expanding.
Recalling that $F \colon Y \to Y$ denotes the induced map, 
it is standard \cite[Cor. p.~199]{AD}, \cite[Proof of Thm.~1]{Young99} that there is
a unique absolutely continuous $F$-invariant probability measure
$\nu_Y$ on $Y$ with $\frac1c \le d\nu_Y / dm\le c$ for some $ c>0$. Let $\nu$ be the corresponding $T$-invariant probability measure  on $X$.

We shall say that the return times of $T$ have a subexponential (or stretched exponential)  moment of order $\gamma \in \,  ]0,1]$, if $\int {\rm e}^{c R^{\gamma}} d m < \infty$ for some $c>0$.
%

The aim of this section is to provide estimates of the $\tau$-mixing coefficients associated with $(\varphi (T^i ) )_{i \geq 0}$  for  a mesurable function $\varphi : (X, \nu) \rightarrow {\mathbb R}$. More precisely, these coefficients are defined by 
\begin{equation} \label{deftaureverse}
  \tau_\varphi (n) =\sup_{\ell \geq 1}  \frac 1 \ell \sup_{1\leq i_1 < \cdots < i_\ell} \tau(\sigma( \varphi(T^{j}), j \geq i_\ell +n), \varphi(T^{i_1}), \ldots, \varphi (T^{i_\ell}))
\end{equation}
where, for any $Z$  taking values in ${\mathbb R}^k$ and any $\sigma$-algebra ${\mathcal M}$, 
\[
\tau ( {\mathcal M}, Z) = \left  \Vert  \sup \Big \{ \Big |  \int f(x) {\mathbb P}_{Z | { \mathcal M} } (dy) -  \int f(x) {\mathbb P}_{Z} (dy) \Big | , f \in \Lambda_1 ({\mathbb R}^k)  \Big \}  \right \Vert_1 \, .
\]
Above $\Lambda_1 ({\mathbb R}^k) $ is the space of $1$-Lipschitz functions from ${\mathbb R}^k$ to ${\mathbb R}$ with respect to the $\ell^1$ distance on ${\mathbb R}^k$. Note that compared to the  $\tau$-mixing coefficients as defined in \cite[Section 7]{DP2005}, there is a time inversion in the definition of  the coefficients defined in \eqref{deftaureverse}.

For an H\"older function $\varphi$ from $X$ to ${\mathbb R}$,  with H\"older's index $\eta$, let
$$
\vert \varphi \vert_{\eta} = \sup_{x,y \in X} \frac{|\varphi(x) - \varphi(y)| }{ d(x,y)^{\eta}} \quad \text{and}  \quad \Vert \varphi  \Vert_{\eta}  = \Vert \varphi  \Vert_{\infty} + \vert \varphi  \vert_{\eta} 
$$
In case where $X=[0,1]$, and $\varphi$ is a bounded variation (BV) function, let also
$$
   \| \varphi \|_v=  \Vert \varphi  \Vert_{\infty} +\|d \varphi \|\, ,
$$
where  $ \|d \varphi \| $ is the variation norm of the measure $d \varphi$.

\begin{Proposition} \label{controltau} Assume that  $T$ is nonuniformly expanding with return time to the basis having a subexponential moment of 
 order $\gamma \in \,  ]0,1]$.  Then, on $(X, \nu)$, one has: 
\begin{enumerate} 
\item  Let $\varphi$ be an H\"older observable with H\"older's index $\eta$, then there exist two positive constants $c_1, c_2$ such that 
\begin{equation} \label{majtauT}
  \tau_\varphi (n)  \leq c_1   \Vert \varphi \Vert_{\eta }e^{-c_2 n^\gamma} \, .
\end{equation}
\item  Let $X=[0,1]$ and $\varphi$ be a BV observable from $[0,1] $ to $[0,1] $.  Assume in addition that there exist $\eta  \in \,  ]0,1]$ and $C>0$ such that, for any $x,y \in [0,1]$ with $x \leq y$, 
\beq  \label{hyposurnu}\nu ([x,y] ) \leq  C (y-x)^\eta \, .
\eeq 
Then there exist two positive constants $c_1, c_2$ such that 
\begin{equation} \label{majtauT2}
  \tau_\varphi (n)  \leq c_1   \Vert \varphi \Vert_v e^{-c_2 n^\gamma} \, .
\end{equation}

\end{enumerate}
\end{Proposition}
\begin{Remark} Following the lines of the proof of Proposition \ref{controltau} and taking into account Lemma 5.7, Proposition 3.2 and inequality (3.6) in \cite{CDKM1}, we infer that if $T$ is nonuniformly expanding with return time to the basis having a  moment of  order $\beta >1$ (meaning that $\int R^{\beta} dm  < \infty$), then, for $\varphi$ an H\"older observable, one has  $\sum_{k \geq 1} k^{\beta -2}   \tau_\varphi (k) < \infty $. 
\end{Remark}
\begin{Remark} Assumption \eqref{hyposurnu} in Item 2 of the proposition above is  satisfied for the class of transformations described  in \cite{Ho05}, with the additional restriction that 
$\lim_{x \rightarrow 0} x^{\varepsilon} \phi(x) = \infty$, for some $\varepsilon >0$. In particular it will be true for the example of Section \ref{sec:example}.
\end{Remark}

\section{Deviations inequalities and Moderate Deviation Principle for Birkhoff sums} \label{DIMDP}

\setcounter{equation}{0}

Starting from  Proposition \ref{controltau} and using the deviations inequality stated in \cite{MPR12,MPR09}, the following corollary holds. 

\begin{Corollary} \label{ineexpo} Assume that  $T$ is nonuniformly expanding with return time to the basis having a subexponential moment of 
 order $\gamma \in \,  ]0,1]$.    Let $\varphi$ be an H\"older observable.
 Let 
$S_n( \varphi) = \sum_{i=0}^{n-1}  \varphi  (T^i) - n \nu (\varphi) $ and 
$V=   {\rm Var}_\nu  (\varphi)+   2 \sum_{i \geq 1}|  {\rm cov}_\nu ( \varphi, \varphi \circ T^i ) |$. Then, on $(X,\nu)$,  \begin{enumerate} 
\item  if $ \gamma \in \,  ]0,1[$ then  for any $n\geq 4$, there exist positive
constants $C_{1}$, $C_{2}$, $C_{3}$ and $C_{4}$ depending on $%
(\gamma, C) $ such that, for any positive $x$,
\begin{multline*} \label{resultineexpo}
\nu \Bigl(\sup_{j\leq n}|S_{j} ( \varphi) |\geq x\Bigr)\leq n\exp \Bigl(-\frac{%
x^{\gamma }}{C_{1}}\Bigr )+\exp \Bigl(-\frac{x^{2}}{C_{2}(1+nV)}\Bigr)  +\exp %
\Bigl(-\frac{x^{2}}{C_{3}n}\exp \Bigr(\frac{x^{\gamma (1-\gamma )}}{%
C_{4}(\log x)^{\gamma }}\Bigr)\Bigr)\,.
\end{multline*}

\item  if $\gamma=1$,  then there exist positive
constants $C_{1}$, and $C_{2}$ such that, for any positive $x$,
\begin{equation*}
\nu \Bigl(\sup_{j\leq n}|S_{j} ( \varphi) |\geq x\Bigr)\leq C_1 \exp \Bigl(-\frac{  x^{2}}{ C_2 (1 +n V  + x  (\log n)^2 ) }\Bigr)\,.
\end{equation*}
\end{enumerate}
In addition, if  $X = [0,1]$,  $\varphi$ is a  BV observable  and  $\nu$ satisfies \eqref{hyposurnu}, then the upper bounds given in Items 1 and 2 are satisfied. 
\end{Corollary}
\begin{Remark} 
Item 2 of Corollary \ref{ineexpo} without the supremum is a consequence of  Inequality (2.3) of \cite[Theorem 2]{MPR09}. The fact that  it also holds  for the maximum $\sup_{j\leq n}|S_{j} ( \varphi) |$ follows from  \cite[Example 2]{KM13}. 
\end{Remark}

If $\gamma \in \,  ]0,1[$, we deduce  from Corollary \ref{ineexpo} that the following large deviation estimate holds: for any $x>0$,
\begin{equation}\label{large-deviation}
\limsup_{n \rightarrow \infty} \frac{1}{n^\gamma}\log \nu  \Big (\frac{S_n( \varphi )}{n}  \geq x  \Big ) \leq \frac{-x^\gamma}{C_1} \, .
\end{equation}
It follows from Section \ref{sec:example} that the exponent of $n$ 
in \eqref{large-deviation} is optimal under our assumption on the return times. 
 \begin{Remark} \label{RemarkonAF}
Note that  an  application of Theorem 2 in  \cite{AF} gives the following large deviation estimate:   for any $x>0$,
\[
\limsup_{n \rightarrow \infty}  \frac{1}{n^{\gamma/(\gamma+1)}} \log \nu \Big (\frac{S_n( \varphi )}{n}  \geq x  \Big ) \leq \frac{-x^{2\gamma/(\gamma+1)}}{C_1}  \, ,
\]
giving then a suboptimal power of $n$. 
    \end{Remark}

\medskip

We also deduce  that the following moderate deviation estimate holds: 
if $a_{n}\rightarrow 0$ and $a_{n}n^{\gamma
/(2-\gamma )}\rightarrow \infty $, then for any positive $x$, 
\beq\label{md1}
\limsup_{n \rightarrow \infty} a_{n}\log\nu \Big (\sqrt{\frac{a_{n}}{n}} S_n( \varphi )  \geq x  \Big ) \leq \frac{-x^2}{C_2 V} \, .
\eeq

In what follows, we shall give a more precise result than \eqref{md1}, and show that $(S_n( \varphi ) / \sqrt n )$ satisfies a Moderate Deviation Principle (MDP). Let us recall the definition.  Let $(Z_{n})_{n}$ be a sequence of random variables defined on $(X, \nu)$. 
We say that the MDP holds for $(Z_{n})_{n}$  
with speed $a_{n}\rightarrow 0$ and good rate function $I(\cdot)$, if the level sets $\{x, I(x) \leq \lambda \}$ are compact for all $\lambda < \infty$, and for each Borel set $A$,
\begin{equation} \label{mdpdef}
-\inf_{t\in A^{o}}I(t) \leq \liminf_{n}a_{n}\log  \nu (\sqrt{a_{n}}%
Z_{n}\in A)  
\leq \limsup_{n}a_{n}\log \nu (\sqrt{a_{n}}Z_{n}\in A)\leq
-\inf_{t\in \bar{A}}I(t)\,,  
\end{equation}
where $\bar{A}$ denotes the closure of $A$ and $A^{o}$ the interior of $A$.

Proposition \ref{controltau} together with Corollary 1 in \cite{MPR12} give the following result.

\begin{Corollary} \label{MDP} Assume that  $T$ is nonuniformly expanding with return time to the basis having a subexponential moment of 
 order $\gamma \in \,  ]0,1[$.    Let $\varphi$ be an H\"older observable. Let 
$S_n( \varphi) = \sum_{i=0}^{n-1}  \varphi  (T^n) - n \nu (f) $ and $\sigma _{n}^{2} =\mathrm{Var}S_{n} (\varphi)$. Assume
 that $\sigma _{n}^{2}\rightarrow \infty $. Then $%
\lim_{n\rightarrow \infty }\sigma _{n}^{2}/n=\sigma ^{2}>0$. Moreover, for all
positive sequences $a_{n}$ with $a_{n}\rightarrow 0$ and $a_{n}n^{\gamma
/(2-\gamma )}\rightarrow \infty $, $\{n^{-1/2}S_{n}(\varphi)\}$ satisfies (\ref%
{mdpdef}) with the good rate function $I(t)=t^{2}/(2\sigma ^{2})$.  
In addition, this Moderate Deviation Principle still holds if  $X = [0,1]$, $\varphi$ is a  BV observable  and  $\nu$ satisfies \eqref{hyposurnu}.
\end{Corollary}
\begin{Remark}
If $\gamma=1$ and $(a_n)_{n \geq 1}$ satisfies $a_{n}\rightarrow 0$ and $n a_{n} \rightarrow \infty $, as $n \rightarrow \infty $, the moderate deviations principle  holds  for H\"older observables, as shown in  \cite[Theorem 4.6]{RY}.  However, as far as we know, this result has not yet been proved when $X = [0,1]$ and $\varphi$ is a  BV observable. But, applying Corollary 5 in \cite{MPR09} together with Proposition \ref{controltau}, one gets (assuming also \eqref{hyposurnu})  that, for all
positive sequences $a_{n}$ with $a_{n}\rightarrow 0$ and $n a_{n} / ((\log n)^2 (\log\log n)^2)\rightarrow \infty $, $\{n^{-1/2}S_{n}(\varphi)\}$ satisfies (\ref%
{mdpdef}) with the good rate function $I(t)=t^{2}/(2\sigma ^{2})$. 
\end{Remark}

\section{Concentration inequalities} \label{Sec:CI}

\setcounter{equation}{0}

In this section, we follow the approach of Chazottes-Gou\"ezel \cite{CG12} to derive concentration inequalities in case of nonuniformly expanding transformations with return time to the basis having a subexponential moment. 

\medskip
Recall that a function $K:X^n \rightarrow {\mathbb R}$ is said to be  a separately H\"older function on $X^n $ of order $\eta \in [0,1]$ if for all $i$ there exists a constant $L_i$ with 
\[
\left | K( x_1, \dots, x_{i-1}, x_i, x_{i+1}, \ldots, x_n )  - K( x_1, \dots, x_{i-1}, x'_i, x_{i+1}, \ldots, x_n ) \right |  \leq L_i d^\eta(x_i, x_i' ) \, , 
\]
for all points $x_1, \ldots, x_n, x'_i$ in $X$. 
\medskip

The following deviation bound holds. 

\begin{Theorem} \label{Concentration} Assume that  $T$ is nonuniformly expanding with return time to the basis having a subexponential moment of 
 order $\gamma \in \,  ]0,1[$.   Let $K$ be a separately H\"older function on $X^n $. Then, there exists a positive constant $\kappa$ such that  for any positive integer $n$ and any $t >0$, 
\begin{equation} \label{ineconc}
\nu  \Big ( x :    K (x, T(x), \dots, T^{n-1} (x) )  -  {\mathbb E} (K )   \geq t   \Big )  \leq 2 \exp  \left (  - \frac{t^2}{\kappa \big  ( \sum_{i=0}^{n-1}L_i^2 +1 + t^{2- \gamma}  \big ) }\right ) 
\end{equation}
where ${\mathbb E} (K )  =    \int K (x, T(x), \dots, T^{n-1} (x)  )  d \nu (x) $. 
\end{Theorem}
\begin{Remark} To prove \eqref{ineconc}, the assumption  \eqref{nonpathologicalH} on $(T,m)$ is not needed. 
\end{Remark}
\begin{Remark} When $\gamma=1$,  Chazottes and Gou\"ezel \cite{CG12} have proved the following concentration inequality:  there exists a positive constant $\kappa$ such that  for any positive integer $n$ and any $t >0$, 
\begin{equation} \label{ineconcCG}
\nu  \Big ( x :    K (x, T(x), \dots, T^{n-1} (x) )  -  {\mathbb E} (K )   \geq t   \Big )  \leq  \exp  \left (  - \frac{t^2}{\kappa  \sum_{i=0}^{n-1}L_i^2  }\right )  \, .
\end{equation}

\end{Remark}

\section{An example of nonuniformly expanding system with stretched exponential return times}
\label{sec:example}

\setcounter{equation}{0}

Let us consider the following examples 
of interval maps, as defined in \cite[Annex A]{CDKM2}, whose return times to the basis satisfy $m(R> n) \sim \mathrm{e}^{-\kappa n^\gamma}$,
where $\gamma \in \,  ]0,1]$ is a parameter and $\kappa = \kappa(\gamma) >0$. 

Let \(T \colon [0,1] \to [0,1]\),
\begin{equation} \label{eq:LSV}
 T(x) = \begin{cases}
    x\bigl(1 + \frac{c}{|\log x|^\beta }\bigr), & x \leq 1/2 \\
    2 x - 1, & x > 1/2
  \end{cases}
\end{equation}
with $\beta = \gamma^{-1} - 1$ and $c=(\log 2)^\beta$ so that $T(1/2) = 1$. Note that this map belongs to the class of transformations considered in Holland \cite{Ho05}. 

Let $Y = ]1/2,1]$ be a base, $R\colon Y \to \N$, $R(x) = \inf \{k \geq 1 : T^k(x) \in Y\}$
be the first return time and $F \colon Y \to Y$, $F(x) = T^{R(x)}(x)$ be the induced map.
Let $\Gamma$ denote the partition of $Y$ into the intervals where $R$ is constant.
Let $m$ denote the Lebesgue measure.

Items a)-d) of the next result have been established in  \cite[Theorem  A.A.]{CDKM2}. To prove its  last part concerning the invariant density, it suffices to follow the lines of the proof of \cite[Lemma 3]{Ho05}. 

\begin{Theorem} 
  \label{thm:ex}
  $T$ is a nonuniformly expanding map with basis $Y$, return time $\tau$ and
  reference measure $m$. That is, there exists $C > 0$ such that
  for every $a \in \Gamma$ and all $x,y \in a$,
  \begin{enumerate}
    \item[a)]\label{thm:ex:bi} $F \colon a \to Y$ is a nonsingular bijection;
    \item[b)]\label{thm:ex:exp} $F$ is expanding: $|F(y) - F(x)| \geq 2 |y-x|$;
    \item[c)]\label{thm:ex:dist} $F$ has bounded distortion: $|\log F'(y) - \log F'(x)| \leq C |F(y) - F(x)|$.
  \end{enumerate}
  Further, there exist $\eta_1, \eta_2 > 0$ such that for all $n \geq 1$,
  \begin{enumerate}
    \item[d)]\label{thm:ex:tails}
    $
      {\rm e}^{-\eta_2n^\gamma }
      \leq m(R \geq n) 
      \leq  {\rm e}^{-\eta_1n^\gamma }
    $.
  \end{enumerate}
    In addition, there exists a unique $T$-invariant measure $\nu$ absolutely continuous with respect to the Lebesgue measure, whose density 
  $\varphi$  is such that: $\varphi$ is  lower bounded by a strictly positive constant on $[0,1]$, $\varphi$ is upper bounded  on any interval $[\varepsilon, 1]$ with $\varepsilon >0$,  and  there exist positive reals $\delta$, $c_1$ and $c_2$ such that, $c_1 \leq  |\log (x) |^{-\beta} \varphi(x)  \leq c_2$ on $]0, \delta]$. 
\end{Theorem}

\begin{Theorem} 
  \label{thm:opt}
  Let $T$ be defined by \eqref{eq:LSV}.  For any measurable function $f : [0,1] \rightarrow {\mathbb R}$, let $S_n(f) = \sum_{k=0}^{n-1}f ( T^k  )  - n \nu (f) $. 
  \begin{enumerate}
    \item[a)]\label{thmoptires1}  There exist a Lipschitz function $f$ and  positive constants $c_1,c_2, \kappa_1$ and $\kappa_2$ such that, for $n > 2/ \nu(f)$, 
    \[
  c_1 {\rm e}^{- \kappa_1 n^{\gamma}} \leq   \nu \big ( |S_n(f)| > n \nu (f) /2 \big )   \leq c_2 {\rm e}^{- \kappa_2 n^{\gamma}}  \, .
    \]
    \item[b)]  There exist an unbounded measurable function $f$ satisfying $ A_1 {\rm e}^{- B_1 t^{1/\delta}}\le \nu ( | f | > t ) \leq  A_2 {\rm e}^{- B_2 t^{1/\delta}} $ for some  $A_1,A_2,B_1,B_2 >0$ and $ \delta  >0$,  and positive constants $c_1,c_2, \kappa_1$ and $\kappa_2$ such that
      \begin{equation}\label{bound-example}
  c_1 {\rm e}^{- \kappa_1 n^{\gamma / (1 + \gamma \delta )}} \leq   \nu \big ( |S_n(f)| > n  \big )   \leq   c_2 {\rm e}^{- \kappa_2 n^{\gamma / (1 + \gamma \delta )}}  \, .
    \end{equation}  
      \end{enumerate}
    \end{Theorem}
    \begin{Remark}  As it will be clear from the proof, for the function of Item b), one can take $f(x) = |\log (x) |^{\delta}$. For this function, following  the proof of Item b) and applying  Theorem 1 in \cite{MPR12}, we derive that  the following deviation Inequality holds: for any $n\geq 4$, there exist positive
constants $C_{1}$, $C_{2}$, $C_{3}$ and $C_{4}$ depending on $%
(\gamma, C) $ such that, for any positive $x$,
\begin{multline*}
\mathbb{P}\Bigl(\sup_{j\leq n}|S_{j} ( f) |\geq x\Bigr)\leq n\exp \Bigl(-\frac{%
x^{\gamma_0 }}{C_{1}}\Bigr )+\exp \Bigl(-\frac{x^{2}}{C_{2}(1+nV)}\Bigr)  +\exp %
\Bigl(-\frac{x^{2}}{C_{3}n}\exp \Bigr(\frac{x^{\gamma_0 (1-\gamma_0 )}}{%
C_{4}(\log x)^{\gamma_0 }}\Bigr)\Bigr)\, ,
\end{multline*}
where $ \gamma_0 =  \gamma  ( 1 + \gamma \delta )^{-1}$ and  $V=   {\rm Var}_\nu  (f)+   2 \sum_{i \geq 1}|  {\rm cov}_\nu ( f,  f \circ T^i ) |$.  From \cite{MPR12}, one can also get a moderate deviations principle for $S_{n} ( f) $. 

Some observables as $ x \mapsto  |\log (x) |^{\delta}$  are considered by Nicol and T\"or\"ok \cite{NT}. In particular, they obtained the lower bound of Item b) for the doubling map which corresponds to $\gamma=1$ 
in \eqref{eq:LSV}.  
    \end{Remark}

\section{Proofs} \label{Sec:proofs}

\setcounter{equation}{0}

\subsection{ Proof of Proposition \ref{controltau}} From Corollary 2.5  in \cite{CDKM1}, recall that 
$$
(\varphi(T^n))_{n \geq 0}= (\psi(g_n, g_{n+1}, \ldots ))_{n \geq 0} \quad \text{in distribution,}
$$
where $(g_i)_{i \geq 0}$ is a strictly stationary Markov chain. This chain is generated by a random variable 
$g_0$ and a sequence of iid innovations $(\varepsilon_i)_{i \geq 1}$ independent of $g_0$. Let  
${\mathcal F}_k= \sigma(g_i, i \leq k)$ and ${\mathcal G}_m= \sigma(g_i, i \geq m)$. The chain is also $\beta$ mixing in the sense that
\begin{equation} \label{But1taup1}
\sup_{k\geq 0} \beta({\mathcal F}_k, {\mathcal G}_{k+n})= O(e^{-Bn^\gamma})\, .
\end{equation}
(See  \cite[Chapter 3]{Br07} for a definition of the $\beta$-mixing coefficients, Relation (3.6) in \cite{CDKM1} and Lemma 2.1 in \cite{CDKM2} for the upper bound \eqref{But1taup1}). 
Let $X_n= \psi(g_n, g_{n+1}, \ldots )$.
We want to prove that 
\begin{equation} \label{But1tau}
  \tau(n)= \sup_{\ell \geq 1}  \frac 1 \ell \sup_{1\leq i_1 < \cdots < i_\ell} \tau({\mathcal G}_{i_\ell +n}, X_{i_1}, \ldots, X_{i_\ell})=O(e^{-Cn^\gamma}) \, ,
\end{equation}
which will imply \eqref{majtauT}. To do so, we need an independent copy $(\varepsilon_i')_{i \geq 1}$ of $(\varepsilon_i)_{i \geq 1}$, this copy being also independent of $g_0$. 
Let now 
$$
 X'_{i_k}= \psi(g_{i_k}, g_{i_k+1}, \ldots, g_{i_k+[n/2]}, g'_{i_k+[n/2]+1}, g'_{i_k+[n/2]+2}, \ldots) \, ,
$$
where the $'$ means that we have used the innovations $(\varepsilon_i')_{i \geq 1}$ to continue the trajectory of the chain. 
From \cite[Proposition 2.3]{CDKM2}, we know that 
$$
\|X'_{i_k}- X_{i_k}  \|_1 \leq K |\varphi|_{\eta}e^{-\delta n^\gamma} \, ,
$$
for some constant $\delta$ depending on $\eta$. 
Hence, from the definition of $\tau(n)$, we infer that 
\begin{equation}\label{But1taup1.5}
  \tau (n) \leq K |\varphi|_{\eta}e^{-\delta n^\gamma}   +  \sup_{\ell \geq 1}  \frac 1 \ell \sup_{1\leq i_1 < \cdots < i_\ell} \tau({\mathcal G}_{i_\ell +n}, X'_{i_1}, \ldots, X'_{i_\ell}) \, .
\end{equation}
Now, since $\|\psi\|_\infty = \|\varphi \|_\infty$, 
\begin{equation} \label{But1taup2}
\sup_{\ell \geq 1}  \frac 1 \ell \sup_{1\leq i_1 < \cdots < i_\ell} \tau({\mathcal G}_{i_\ell +n}, X'_{i_1}, \ldots, X'_{i_\ell}) \leq  \|\varphi \|_\infty   \sup_{k\geq 0} \beta({\mathcal F}_k', {\mathcal G}_{k+[n/2]})
\end{equation}
where ${\mathcal F}_k'={\mathcal F}_k \vee \sigma(\varepsilon_i', i \geq 1)$, and we have used that $\beta({\mathcal F}_k', {\mathcal G}_{k+[n/2]})=\beta( {\mathcal G}_{k+[n/2]}, {\mathcal F}_k')$. Since $(\varepsilon_i')_{i \geq 1}$ is independent of $(g_i)_{i \geq 0}$, standard arguments (see \cite[Theorem 6.2]{Br07}) show that
\begin{equation} \label{But1taup3}
\beta({\mathcal F}_k', {\mathcal G}_{k+[n/2]}) \leq \beta({\mathcal F}_k, {\mathcal G}_{k+[n/2]}) \, .
\end{equation}
The estimate \eqref{But1tau}  follows from \ \eqref{But1taup1},  \eqref{But1taup1.5}, \eqref{But1taup2} and \eqref{But1taup3}.

\medskip

We turn now to Item 2; namely, the case of a BV function $\varphi$ on $[0,1]$, when $T: [0,1] \rightarrow [0,1]$ is a nonuniformly expanding map  with return time to the basis satisfying 
$m(R>n) = O(e^{-\delta n^\gamma})$ for $\gamma \in  \, ]0,1]$ and some $\delta >0$. We want to control $\tau_\varphi (n)$,
and we recall that 
$$
\varphi(x)=  \varphi(0)+\int_0^x d\varphi(t) \, .
$$
Let now 
$$
f_{\epsilon}(t,x)= {\bf 1}_{t \leq x} + \frac 1 \epsilon (x + \epsilon -t){\bf 1}_{ x<t\leq x + \epsilon}\
\quad \text{and} \quad  \varphi_\epsilon(x)= \varphi(0)+\int_0^1 f_{\epsilon}(t,x) d\varphi(t)\, .
$$
Since $\nu$ satisfies \eqref{hyposurnu}, note that
$$
   \|\varphi-\varphi_\epsilon \|_{1, \nu} \leq C_\eta \|d \varphi \|  \epsilon^\eta \, .
$$
On the other hand, using that $u\mapsto f_\varepsilon (t,u)$
is $\frac1\varepsilon$-Lipschitz, we have
$$
  |\varphi_\epsilon (x) - \varphi_\epsilon (y) | \leq \int_0^1 |f_\epsilon(t,x)-f_\epsilon(t,y)| |d\varphi| (t) \leq 
   \frac{\|d\varphi\|}{\epsilon} |x-y| \, .
$$
From the definition of $\tau_\varphi(n)$ and Item 1, we infer that 
\begin{equation*}
  \tau_\varphi  (n) \leq  C_\eta \| \varphi \|_v \epsilon^\eta +  \frac{K \|\varphi\|_v  (1 + \epsilon)  }{\epsilon}  e^{-Cn^\gamma} \, .
\end{equation*}
The upper bound  \eqref{majtauT2} easily follows, with
$c_2= (C \eta)/(\eta +1)$.

\subsection{Proof of Theorem \ref{Concentration}}   First, we recall that we can associate to the transformation  $T$ a Young tower ${\mathcal T}$ and a transformation ${\bar T}$. More precisely, ${\mathcal T}$ is the space
\[
{\mathcal T}=\{(y,i): y \in Y, i < R(y)\}
\]
 and the map $\bar T$ on ${\mathcal T}$ is defined by 
\[
 \bar T(y,i)=\begin{cases}
  (y, i+1) \quad \quad \quad   \text{if $i< R(y)-1 $}\\
  (T^{R(y)}(y), 0) \quad  \text{if $i= R(y)-1$.}
  \end{cases}
\]
For any $\alpha \in E$, define the height $h_{\alpha}$ by  $h_{\alpha} = R(y)$ for $y \in \Gamma_{\alpha}$. One can then  define the  floors of the tower $\Delta_{\alpha,i}$ for
 $ \alpha \in E$ and $i \in \{0, \ldots, h_{\alpha} -1\}$:
 $\Delta_{\alpha,i}= \{ (y, i): y \in  \Gamma_\alpha \}$.
 These floors define a partition of ${\mathcal T}$:
\[
{\mathcal T}= \bigcup_{ \alpha \in E, i \in \{0, \dots , h_{\alpha}-1\}} \Delta_{\alpha,i} \, .
\]
On the tower, there exists a reference measure $\bar m$ defined as follows:  if $\bar B$ is a set included in $\Delta_{\alpha,i}$, that can be written
 as $\bar B= B\times \{ i \}$ with $B \subset \Gamma_{\alpha}$,
 then ${\bar m}(\bar B)=m(B)$.
From Young \cite{Young99} (see also \cite{Go02}, Proposition 1.3.18)), it follows that   on the tower, there exists a unique $\bar T$-invariant probability measure $\bar \nu$
 which is absolutely continuous with respect to $\bar m$.  On another hand, the distance on the tower is defined by $ \delta(x,y)=\lambda^{-s(x,y)}$ where $s(x,y)$ is the separation time, i.e. the number of returns to the basis   $\bar Y=\{(y,0), y \in Y\} $ before the iterates of the points $x$ and $y$ are not in the same element of the partition.
Let now $\pi$ be the ``projection'' from ${\mathcal T}$ to $X$ defined by
 $
 \pi(y,i)=T^i(y)
 $. 
 Then, one has $ \pi \circ \bar T= T \circ \pi$ and for any $x, y $ in ${\mathcal T}$, there exists $C>0$, such that 
\[
d( \pi(x) , \pi(y) ) \leq C \delta (x,y)   \, .
\]

Moreover it can be checked that  the $T$-invariant measure $\nu$ defined in Section \ref{Sec:mixing},  is the image measure of $\bar \nu $ by $\pi$.   Let ${\bar R}$ be the function from $\bar Y$ to $R(Y)$ such that ${\bar R} (y,0)= R(y)$. The quantity $\bar \nu (\{(y, 0) \in \bar Y :
  {\bar R}((y,0))>k\})$ is exactly of the same order as $
   m (\{y \in Y : R(y)>k\})$ (see \cite{Go02}, Proposition 1.1.24). From all these considerations, we infer that if $T$ is nonuniformly expanding then ${\bar T}$ is also nonuniformly expanding with respect to the distance  $ \delta^\eta$ whatever $\eta \in \,  ]0,1]$.  On another hand, recalling that $\nu$ is the image measure of $\bar \nu $ by $\pi$, and that $ \pi \circ \bar T= T \circ \pi$, we get, for any $t >0$, 
\begin{multline*}
\nu  \Big ( x :   K (x, T(x), \dots, T^{n-1} (x) )  -  {\mathbb E} (K )   \geq t   \Big )  \\
=  {\bar \nu }   \Big ( z :   K ( \pi(z), T \circ \pi (z), \dots, T^{n-1} \circ \pi (z) )  -  {\mathbb E} (K )   \geq t   \Big ) \\
=  {\bar \nu }   \Big ( z :   K ( \pi(z),  \pi  \circ {\bar T}(z), \dots,  \pi \circ {\bar T}^{n-1} (x) )  -  {\mathbb E} (K )   \geq t   \Big )  \, .
\end{multline*}
Next, defining ${\tilde K} $ on ${\mathcal T}^{\mathbb N}$ by 
\[
{\bar  K} (z_1, \ldots, z_n) = K ( \pi(z_1), \ldots, \pi(z_n) ) \, , 
\]
one sees that it satisfies
\[
 \Big | {\bar K} (z_1, \ldots, z_n)  -  {\bar K} (z'_1, \ldots, z'_n) \Big | \leq \sum_{i=1}^n L_i \delta^\eta (z_i,z_i') \, . 
\]
So ${\bar  K}$ is separetely Lipschitz with respect to $\delta^\eta$.  Hence, it suffices to proof Theorem \ref{Concentration} with ${\bar  K}$ instead of $K$,  ${\bar T}$ instead of $T$ and  ${\bar \nu } $ instead of $\nu$. Since the tower is also nonuniformly expanding with respect to $\delta^\eta$, we shall only consider the case $\eta=1$   as in \cite{CG12}. 
%

From now on, we shall use the same notations as in  \cite{CG12,GM14}. Hence, we shall now use $T,K,\nu$ instead of  ${\bar T},{\bar K},{\bar \nu}$

It is convenient to consider $K$ as a function defined on the space ${\tilde {\mathcal T}}= {\mathcal T}^{{\mathbb N}}$ endowed with the probability measure 
${\tilde \nu} = \nu \otimes \delta_{Tx} \otimes \delta_{T^2x} \otimes \cdots$. Let ${\mathcal F}_k $ be the $\sigma$-algebra generated by indices starting with $k$. Let 
\[
K_k (x_k, \ldots )= \E ( K | {\mathcal F}_k ) ( x_k, \ldots)  = \sum_{T^k x = x_k} g^{(k)} (x)  K 
\]
where $g^{(k)}  (x)  = g(x)\cdots g(T^{k-1} (x))$ is the inverse of the jacobian of $T^k$. Next define $D_k = K_k - K_{k+1}$. Clearly $D_k $ is ${\mathcal F}_k $-measurable and such that 
$\E ( D_k | {\mathcal F}_{k+1})=0$. Hence $(D_k)_{k \geq 1}$ is a sequence of reversed martingale differences w.r.t. the decreasing filtration $({\mathcal F}_k)_{k \geq 1} $.

As quoted after the statement of \cite[Theorem 3.1]{CG12}, we can assume without loss of generality  that $\sup_{i \geq 0} L_i \leq \varepsilon_0$ for some $\varepsilon_0 >0$ (the appropriate $\varepsilon_0$ will be chosen later and will appear in the constant $\kappa$ of inequality \eqref{ineconc}). 
Noticing that $K - \E (K) = \sum_{k=0}^n D_k $, inequality \eqref{ineconc} will follow from  the reversed martingale differences sequences version  of Theorem 2.1 in Fan et al  \cite{FGL17} (see also their inequality (1.5)).  Note that \cite[Theorem 2.1]{FGL17} is stated for sums of martingale differences but with the same proof it also holds for partial sums associated with reversed martingale differences sequences. Let us state it for reader convenience. Assume that $(D_p)_{p \geq 1}$ is  a sequence of reversed martingale  differences  with  respect  to  the  filtration $({\mathcal F}_p)_{p \geq 1}$ (so $D_p$ is ${\mathcal F}_p$-adapted and such that ${\mathbb E} (D_p | {\mathcal F}_{p+1})=0 $ a.s.), then setting 
\[
u_n = \max \Big (  \Big \Vert  \sum_{p=1}^n \E \big ( D_p^2 {\rm e}^{(D_p^+)^{\gamma}}| {\mathcal F}_{p+1}  \big )  \Big \Vert_{\infty}, 1 \Big ) \, , 
\]
one has, for any $x>0$,
\[
{\mathbb P} \Big ( \max_{1 \leq k \leq n} \sum_{i=1}^kD_i \geq  x \Big )  \leq 2  \exp \Big \{  - \frac{x^2}{2 (u_n + x^{2 - \gamma } ) } \Big \} \, .
\]
 Hence we need to prove that 
\begin{equation} \label{toverifyconditionFan}
\Big \Vert  \sum_{p=1}^n \E \big ( D_p^2 {\rm e}^{(D_p^+)^{\gamma}}| {\mathcal F}_{p+1}  \big )  \Big \Vert_{\infty} \leq C \sum_{i=1}^n L_i^2  \, .
\end{equation}
As noticed at the beginning of the proof of \cite[Lemma 3.3]{CG12}, if $x_{p+1} \notin \Delta_0= \cup_{\alpha \in E} \Delta_{\alpha, 0}$ (hence when $x_{p+1} $ is not in the basis of the tower), $ \E \big [ D_p^2 {\rm e}^{(D_p^+)^{\gamma}}| {\mathcal F}_{p+1}  \big ] (x_{p+1}, \ldots ) =0$. Now when $x_{p+1} \in  \Delta_0$, denoting by $z_{\alpha}$ its preimages with respective heights $h_{\alpha}$, we have
\[
\E \big [ D_p^2 {\rm e}^{(D_p^+)^{\gamma}}| {\mathcal F}_{p+1}  \big ] (x_{p+1}, \ldots ) = \sum_{\alpha} g(z_{\alpha})  {\rm e}^{(A_p^+(z_{\alpha}) )^{\gamma}} A_p^2(z_{\alpha}) \, , 
\]
where 
\[A_p(z) :=D_p  (z, x_{p+1}, \ldots )  = K_p  (z, x_{p+1}, \ldots ) - K_{p+1}  (x_{p+1}, \ldots ) \, .
\]
According to Lemma 4.2 in \cite{GM14} and taking into account that $m(R>n) = O(e^{-\delta n^\gamma})$, for some $\delta >0$,  it follows that if $z$ is at height $h$ and $x_{k+1} = Tz  \in \Delta_0$, then there exist positive constants $C$ and $\kappa$ such that 

\begin{equation} \label{estimateAz}
|A_p(z)| \leq C \sum_{a=0}^{p-h} L_a  {\rm e}^{- \kappa (p-h -a)^{\gamma}} +  \sum_{a=p-h+1}^{p} L_a \, .
\end{equation}
Since, it is assumed that $\sup_{i \geq 0} L_i \leq \varepsilon_0$,  \eqref{estimateAz} implies that $|A_p(z)| \leq c_0 (h+1) \varepsilon_0$. 
Using again the upper bound \eqref{estimateAz}, it follows that 
\[
 \big  \Vert \E \big [ D_p^2 {\rm e}^{(D_p^+)^{\gamma}}| {\mathcal F}_{p+1}  \big ]  \big \Vert_{\infty}  \ll  \sum_{h \geq 0} \nu (R= h+1)  {\rm e}^{ ( c_0  \varepsilon_0 h )^{\gamma}}\Big (  \sum_{a=0}^{p-h} L^2_a  {\rm e}^{- \kappa (p-h -a)^{\gamma}} +(h+1)  \sum_{a=p-h+1}^{p} L^2_a
  \Big )  \, .
\]
Hence, for $\varepsilon_0 $ such that $( c_0  \varepsilon_0  )^{\gamma}= A/2$, using that $m(R>n) = O(e^{-\delta n^\gamma})$, there exists a positive constant $c_1$ such that 
\begin{multline*}
\big  \Vert \E \big [ D_p^2 {\rm e}^{(D_p^+)^{\gamma}}| {\mathcal F}_{p+1}  \big ] \big  \Vert_{\infty}  \ll  \sum_{h \geq 0}  {\rm e}^{ - c_1  h^{\gamma}}\Big (  \sum_{a=0}^{p-h} L^2_a  {\rm e}^{- \kappa (p-h -a)^{\gamma}} +(h+1)  \sum_{a=p-h+1}^{p} L^2_a   \Big )  \\
 \ll      \sum_{a=0}^{p} L^2_a   \Big \{ \sum_{0 \leq h \leq p-a}  {\rm e}^{ - c_1  h^{\gamma}}  {\rm e}^{- \kappa (p-h -a)^{\gamma}} + \sum_{ h \geq p-a} (h   +1)  {\rm e}^{ - c_1  h^{\gamma}} \Big \} \, .
\end{multline*}
By splitting $ \sum_{0 \leq h \leq p-a} $ into two sums: $ \sum_{0 \leq h \leq [(p-a)/2]} $ and $ \sum_{[(p-a)/2] < h \leq p-a} $, we derive that there exists a positive constant $c_2$ such that 
\[
\big  \Vert \E \big [ D_p^2 {\rm e}^{(D_p^+)^{\gamma}}| {\mathcal F}_{p+1}  \big ] \big  \Vert_{\infty}  \ll   \sum_{a=0}^{p} L^2_a   {\rm e}^{ - c_2  (p-a)^{\gamma}} \, .
\]
Hence
\[
 \sum_{p=1}^n  \big  \Vert \E \big [ D_p^2 {\rm e}^{(D_p^+)^{\gamma}}| {\mathcal F}_{p+1}  \big ] \big  \Vert_{\infty}  \ll  \sum_{a=0}^{n} L^2_a  \sum_{p \geq a}  {\rm e}^{ - c_2  (p-a)^{\gamma}} 
\]
which leads to \eqref{toverifyconditionFan} and ends the proof of the theorem. \hfill $\square$

\subsection{Proof of Theorem \ref{thm:opt}}

We start with some preliminary considerations. 
Let $S:=T_{]0,1/2]}^{-1}$ be the inverse left branch of $T$ and 
 $U:=T_{]1/2,1]}^{-1}$. 
 
 \medskip
 
 For every $n\in \N_0$, set $y_n= S^n (1/2)$, $I_n:= (0,y_n]$ and 
 $J_n :=U (I_n)=(1/2, y_n/2+1/2]$. 
 
 \medskip
 
 It follows from \cite[Annex A]{CDKM2} that there exist
 $\upsilon_1, \, \upsilon_2>0$ such that 
 $$
{\rm e}^{-\upsilon_1 n^\gamma}\le   y_n \le {\rm e}^{-\upsilon_2 n^\gamma} \, .
 $$
 Then, using the control of the density $\varphi$, we infer that there exist $K,\, L>0$ and $\upsilon_3, \, \upsilon_4>0$  such that 
 $$
L {\rm e}^{-\upsilon_3 n^\gamma} \le \nu(I_n)\le K {\rm e}^{-\upsilon_4 n^\gamma}
 $$
 

\medskip

\noindent{\bf Proof of Item a).} The upper bound given in Item a) comes from an application of inequality \eqref{ineconc} when $\gamma > 1$ and of inequality \eqref{ineconcCG} when $\gamma=1$.  To construct a  
Lipschitz function for which the lower bound holds the idea is to take 
the function ${\bf 1}_{(1/2,1]}$ and to extend it to a Lipschitz function. 

For every $x\in (1/2,1]$, set $f(x)=1$ and for every $x\in I_1$ set $f(x)=0$. 
For every $x\in (y_1, 1/2]$ set $f(x)=\frac{x-y_1}{1/2-y_1}$. 

\medskip

Let $n\in \N$ and  $x\in J_n\subset (1/2,1]$. Then, $f(x)=1$ and for every $k\in [\![1,n-1]\!]$, $f\circ T^k(x)=0$. So $S_n(f)(x)=1-n\nu(f)$ 
and for $n>2/\nu(f)$, $|S_n(f)(x)|\ge n\nu(f)/2$. 
So, there exists a positive constant $c>0$ such that 
$$
\nu(|S_n(f)|\ge n\nu(f)/2)\ge \nu(J_n) \geq c y_n \ge c {\rm e}^{-\upsilon_1 n^\gamma}
\, ,
$$
proving the lower bound.

\medskip

\noindent{\bf Proof of Item b).}  Let $\delta>0$ and $f:= ( |\log x|)^\delta$.  Note that 
\[
\nu ( | f | > t ) = \nu ( ]0 ,  {\rm e}^{-t^{1/\delta}} [)  \leq  C \int_0^{{\rm e}^{-t^{1/\delta}}} ( 1 + |  \log  x |^{\beta} ) dx \leq  A  {\rm e}^{- B t^{1/\delta}}  \, .
\]
The lower bound may be proved similarly. 

\medskip
Let us  prove the lower bound in \eqref{bound-example}. We start by noticing that 
$$\nu(f)\le (  \log 2)^\delta\nu((1/2,1])+\sum_{n\in \N_0} 
( \upsilon_1 (n+1)^{\gamma } )^\delta \nu(I_n\backslash I_{n+1})
\le C+ K \sum_{n\in \N} 
(\upsilon_1 (n+1)^{\gamma } )^{\delta} {\rm e}^{-\upsilon_4 n^\gamma}
<\infty\,.
$$
 Let $1\le m\le n$ be integers. 
Let $x\in I_{m}$ (then $T^kx\in I_{m-k}$ for every integer $k\in [0,m-1]$). We have 
\begin{multline*}
S_n(f):=\sum_{k=0}^{n-1}(f\circ T^k (x)-\nu(f))=
\Big(\sum_{k=0}^{n-1} f\circ T^k(x)\Big) -n\nu(f)\\
\ge \sum_{k=0}^{m-1}(\upsilon_2 (m-k)^\gamma)^\delta -n\nu(f)\ge 
\frac{\upsilon_2^\delta m^{1+\gamma\delta}}{1+\gamma \delta} -n\nu(f)\, .
\end{multline*}
Let $\varepsilon>0$ and take \[
m= \Big[\Big(
\frac{1+\gamma\delta}{\upsilon_2^\delta} (\nu(f)+\varepsilon)n\Big)^{1/(1+\gamma\delta)}\Big] \, , \] which is smaller than $n$, for $n$ large enough since $\gamma\delta>0$. It follows that, for every  $x\in I_m$, $S_n (f) (x)\ge \varepsilon n$. Moreover, 
$$
\nu(I_m)\ge L {\rm e}^{-\upsilon_3 m^\gamma} 
\ge C {\rm e}^{-c n^{\gamma/(1+\gamma\delta)}}\, ,
$$
which proves the lower bound.

\smallskip

The upper bound in \eqref{bound-example} will  follow by applying Theorem 1 in \cite{MPR12} with $\gamma_1 = \gamma$ and $\gamma_2 = 1/ \delta$ provided one can prove that $  \tau_{f} (n) =O(e^{-Cn^\gamma}) $ where we recall that $f:= |\log x|^\delta$. With this aim, let $\varepsilon  \in [0,1[$ and define $f_{\varepsilon} (x) = (\log (1/ \varepsilon ))^{\delta} {\bf 1}_{[0, \varepsilon]} (x) + f (x) 
 {\bf 1}_{]\varepsilon, 1]} (x) $.  Assume first that $\delta \geq 1$. In this case, note that $f_{\varepsilon}$ is a Lipshitz function with Lipshitz  constant  $ \varepsilon^{-1}  \delta | \log \varepsilon|^{\delta -1}$. Hence, by Item 1 of Proposition \ref{controltau}, $  \tau_{f_{\varepsilon}} (n) \leq K    \varepsilon^{-1} | \log \varepsilon|^{\delta -1}  e^{-Cn^\gamma} $. On another hand, for $\varepsilon $ small enough,
 \[
 \nu ( | f- f_{\varepsilon} |)  \leq C \int_0^{\varepsilon} (  (\log (1/ x ))^{\delta} -  (\log (1/ \varepsilon ))^{\delta} ) dx \leq {\tilde C} \varepsilon \, .
 \]
 Hence, for $n$ large enough,  choosing  $\varepsilon = e^{-C n^{\gamma}/3}$, and using the definition of $\tau_f$, we get that, if $\delta \geq 1$, 
 \[
 \tau_f (n) \leq {\tilde C} \varepsilon + K    \varepsilon^{-1} | \log \varepsilon|^{\delta -1}  e^{-Cn^\gamma}=  O(e^{-\kappa n^\gamma})  \, .
 \]
Assume now that $\delta \in ]0,1[$. In this case,  since the fonction $x \mapsto |x|^{\delta}$ is $\delta$-H\"older, we derive that $f_{\varepsilon}$ is a $\delta$-H\"older function with H\"older  constant equals to $ \varepsilon^{-\delta}$, meaning that $| f_{\varepsilon} (x) - f_{\varepsilon} (y)  | \leq \varepsilon^{-\delta} |x-y|^{\delta}$. Using again  Item 1 of Proposition \ref{controltau} and arguing as above, the desired upper bound for  $ \tau_f (n) $ follows.



\begin{thebibliography}{9}

\small


\bibitem{AD} Aaronson, J. and Denker, M. Local limit theorems for partial sums of stationary sequences
generated by Gibbs-Markov maps, \textit{Stoch. Dyn.}  {\bf 1}, no. 2 (2001), 193--237.

\bibitem{AF} Aimino, R. and Freitas J. M. Large deviations for dynamical systems with stretched exponential decay of correlations, \textit{Portugaliae Mathematica} {\bf 76} (2019), no 2, 143--152.

\bibitem{AFLV11} Alves, F. J.;  Freitas, J. M.; Luzzatto, S. and  Vaienti, S.  From
rates of mixing to recurrence times via large deviations. \textit{Adv. Math.} {\bf 228} (2011),
no. 2, 1203--1236.

\bibitem{Br07} Bradley, R.C. \textit{Introduction to strong mixing conditions.} Vol. 1. Kendrick Press, Heber City, UT, 2007. xviii+539 pp.

\bibitem{CG12}   Chazottes, J.R. and Gou\"ezel, S. Optimal concentration inequalities for dynamical systems. \textit{Comm. Math. Phys.} \textbf{316} (2012), no. 3, 843--889.

\bibitem{Ch99} Chernov, N. Decay of correlations and dispersing billiards. \textit{J. Statist. Phys.} {\bf 94} (1999) 513--556.

\bibitem{CDKM1}  Cuny, C.; Dedecker, J.; Korepanov, A. and  Merlev\`ede, F. Rates in almost sure invariance principle for slowly mixing dynamical systems. \textit{Ergodic Theory Dynam. Systems} \textbf{ 40} (2020), no. 9, 2317--2348. 

\bibitem{CDKM2} Cuny, C.; Dedecker, J.; Korepanov, A. and  Merlev\`ede, F. Rates in almost sure invariance principle for quickly mixing dynamical systems. \textit{Stoch. Dyn.} {\bf  20} (2020), no. 1, 2050002, 28 pp.

\bibitem{DGM}  Dedecker, J.; Gou\"ezel, S. and Merlev\`ede, F. Large and moderate deviations for bounded functions of slowly mixing Markov chains. \textit{Stoch. Dyn.} {\bf 18} (2018), no. 2, 1850017, 38 pp.

\bibitem{DP2005}  Dedecker, J. and Prieur, C.  New dependence coefficients. Examples and applications to statistics. \textit{Probab. Theory Related Fields} {\bf 132} (2005), no. 2, 203--236.

\bibitem{FGL17} Fan, X.; Grama, I. and  Liu, Q. 
Deviation inequalities for martingales with applications. \textit{J. Math. Anal. Appl.} {\bf 448} (2017), no. 1, 538--566. 

\bibitem{Go02} Gou\"ezel, S. {\em Vitesse de d\'ecorr\'elation et th\'eor\`emes limites pour les applications non uniform\'ement dilatantes}. PhD Thesis (2004).

\bibitem{Gou06} Gou\"ezel, S. Decay of correlations for nonuniformly expanding systems. 
\textit{Bull. Soc. Math. France} {\bf 134} (2006), no. 1, 1--31.

\bibitem{GM14}   Gou\"ezel, S. and  Melbourne, I. Moment bounds and concentration inequalities for slowly mixing dynamical systems. \textit{Electron. J. Probab.} \textbf{19} (2014), no. 93, 30 pp.

\bibitem{Ho05} Holland, M. Slowly mixing systems and intermittency maps. \textit{Ergodic Theory Dynam. Systems} \textbf{25} (2005), no. 1, 133--159.

\bibitem{KM13}  Kevei, P. and  Mason, D.M. A more general maximal Bernstein-type inequality. \textit{High dimensional probability VI,} 55--62, Progr. Probab., 66, Birkh\"auser/Springer, Basel, 2013.



\bibitem{Mel09} Melbourne, I. Large and moderate deviations for slowly mixing dynamical
systems. \textit{ Proc. Amer. Math. Soc.} {\bf 137} (2009), no. 5, 1735--1741.

\bibitem{MN08}  Melbourne, I. and Nicol, M.  Large deviations for nonuniformly hyperbolic systems. \textit{Trans. Amer. Math. Soc.} {\bf 360} (2008), no. 12, 6661--6676.

\bibitem{MPR09}  Merlev\`ede, F.; Peligrad, M. and Rio, E. Bernstein inequality and moderate deviations under strong mixing conditions. \textit{High dimensional probability V: the Luminy volume}, 273–292, Inst. Math. Stat. (IMS) Collect., 5, Inst. Math. Statist., Beachwood, OH, 2009.

\bibitem{MPR12}  Merlev\`ede, F.; Peligrad, M. and Rio, E. A Bernstein type inequality and moderate deviations for weakly dependent sequences. \textit{Probab. Theory Related Fields} \textbf{151} (2011), no. 3-4, 435--474.

\bibitem{NT} Nicol, M. and T\"or\"ok, A. A note on large deviations for unbounded observables.  \textit{Stoch. Dyn.}  \textbf{20} (2020), no. 5, 2050030, 21 pp.

\bibitem{RY} Rey-Bellet, L. and Young, L.S.
Large deviations in non-uniformly hyperbolic dynamical systems. \textit{Ergodic Theory Dynam. Systems} \textbf{ 28} (2008), no. 2, 587–612. 

\bibitem{Via97}Viana, M. Multidimensional nonhyperbolic attractors. \textit{Inst. Hautes \'Etudes
Sci. Publ. Math.} (1997), no. 85, 63--96.


\bibitem{Young99} Young, L.S. Recurrence times and rates of mixing.  \textit{Israel J. Math.} \textbf{ 110}  (1999), 153-188.



\end{thebibliography}
\end{document}